\documentclass[12pt]{amsart}   

\newtheorem{theorem}{Theorem}[section]

\newtheorem{proposition}[theorem]{Proposition}

\newtheorem{example}[theorem]{Example}
\newtheorem{conjecture}[theorem]{Conjecture}
\newtheorem{corollary}[theorem]{Corollary}

\theoremstyle{definition}

\theoremstyle{remark}
\newtheorem{remark}[theorem]{Remark}
\numberwithin{equation}{section}

\newcommand{\ra}{\rightarrow}

\newcommand{\nats}{{\rm I\kern-.1667em{}N}}
\newcommand{\ints}{\mathbb{Z}}
\newcommand{\reals}{{\rm I\kern-.1667em{}R}}

\begin{document}

\vspace{0.5in}

\title[Generalized Balanced Pair Algorithm]%
{A Generalized Balanced Pair Algorithm}

\author{Brian F. Martensen}
\address{The University of Texas at Austin, Department of Mathematics, 1 
University Station/C1200, Austin, Texas 78712}
\email{martense@math.utexas.edu}



\begin{abstract} 
We present here a more general version of the balanced pair
algorithm.  This version works in the reducible case and terminates more
often than the standard algorithm.  We present examples to illustrate this
point.  Lastly, we discuss the features which lead to balanced pair
algorithms not terminating and state several conjectures.
\end{abstract}

\maketitle

\section{Introduction}

The balanced pair algorithm was first introduced by Livshits \cite{liv1, liv2} for checking
pure discrete spectrum for the $\ints$-action of a substitution.  A version of this algorithm
was presented by Sirvent and Solomyak \cite{sirsol} for irreducible substitutions.  For
substitutions of Pisot type, Sirvent and Solomyak also give an explicit relationship between
this algorithm and an overlap algorithm used in \cite{sol} for checking the pure discrete
spectrum of the $\reals$-action on the substitution tiling space.  Recent results of Clark
and Sadun \cite{clasad} give conditions for the conjugacy of the $\reals$-actions on
substitution tiling spaces which otherwise do not differ combinatorially or topologically.  
Their results give an immediate relation between the $\ints$-action on the sequence space and
the $\reals$-action on the tiling space.  The conditions given in \cite{clasad} also allow us 
a way to generalize the balanced pair algorithm to reducible substitutions.  The procedure 
for doing this is described in Section \ref{balpairalgo} below.

The purpose of extending this algorithm two-fold.  First, it creates a way in which the algorithm
will terminate more often.  When the balanced pair algorithm terminates has been much studied.  
Hollander \cite{hol} (see also \cite{holsol}) has shown that the balanced pair algorithm will
terminate for all $2$-symbols Pisot substitutions.  This fact along with the solution to the
coincidence conjecture for $2$-symbols \cite{bardia2}, has shown that all $2$-symbol Pisot
substitutions have pure discrete spectrum.  This has not yet been shown for an arbitrary number of
symbols.  In fact, it is not yet known whether the balanced pair algorithm terminates for all
Pisot substitutions.

The second reason for extending the algorithm to the reducible case is related to collaring 
or rewriting substitutions to obtain new (yet conjugate) substitutions.  Collaring or 
rewriting procedures generally increase the number of symbols which turns irreducible 
substitutions into reducible ones.  It would be beneficial to know that such procedures did 
not change the potential for a balanced pair algorithm to terminate.  In particular, there 
are rewriting procedures which automatically produce coincidences.  Thus the question of pure 
point spectrum (and even the coincidence conjecture) relies entirely on whether these 
reducible systems terminate.  We explore these and other questions in Section \ref{opensec}.

\section{Preliminaries}

Let $\mathcal{A} = \{1, 2, \ldots , n\}$ be a finite alphabet;  $\mathcal{A}^*$ will denote
the collection of finite nonempty words with letters in $\mathcal{A}$.  A {\em substitution}
is a map $\varphi : \mathcal{A}\ra\mathcal{A}^*$ and extends naturally to $\varphi :  
\mathcal{A}^*\ra\mathcal{A}^*$ and $\varphi : \mathcal{A}^\ints\ra\mathcal{A}^\ints$ by
concatenation.

We associate to each word $w$ a population vector $${\bf p}(w)= (p_1(w), p_2(w), \ldots ,
p_n(w))$$ which assigns to each $p_i(w)$ the number of appearances of the letter $i$ in the
word $w$.  To a substitution $\varphi$ there is an associated transition matrix $A_\varphi =
(a_{ij})_{i\in \mathcal{A}, j \in \mathcal{A}}$ in which $a_{ij}=p_i(\varphi(j))$.  Note that
$A_\varphi({\bf p}(w))={\bf p}(\varphi(w))$.  We say a substitution $\varphi$ is {\it
primitive} if $\varphi^m(i)$ contains $j$ for all $i, j \in \mathcal{A}$ and sufficiently
large $m$.  Equivalently, $\varphi$ is primitive if and only if the matrix $A_\varphi$ is
eventually positive (there exists $m$ such that the entries of $A_\varphi^m$ are strictly
positive), in which case $A_\varphi$ has an eigenvalue $\lambda_\varphi$ larger in modulus
than its remaining eigenvalues called the Perron-Frobenius eigenvalue of $A_\varphi$ (and
$\varphi$).

We form a space $\Omega_\varphi$ as the set of {\em allowed bi-infinite words} for $\varphi$.
That is, ${\bf u}\in\Omega_\varphi$ if and only if for each finite subword $w$ of ${\bf u}$,
there are $i\in\mathcal{A}$ and $n\in\nats$ such that $w$ is a subword of $\varphi^n(i)$.  We
give $\Omega_\varphi$ the subspace topology (of $A^\ints$ with the product topology) and
denote the natural shift homeomorphism by $\sigma$.  Then $(\Omega_\varphi, \sigma)$ is a
topological dynamical system which is minimal and uniquely ergodic provided the substitution
is primitive.  We will often use the fact that, due to unique ergodicity, any allowable word
$w$ appears in any ${\bf u}\in\Omega_\varphi$ with a well-defined and bounded positive
frequency.

The spectral type of the measure-preserving transformation $(\Omega_\varphi, \sigma, \mu)$ 
is, by definition, the spectral type of the unitary operator $U_\varphi:f(\cdot )\mapsto 
f(\sigma\cdot )$ on $L^2(\Omega_\varphi, \mu)$.  We say $(\Omega_\varphi, \sigma, \mu)$ has 
{\it pure discrete spectrum} if and only if there is a basis for $L^2(\Omega_\varphi, \mu)$ 
consisting of eigenfunctions for $U_\varphi$.

Additionally, to each substitution $\varphi$ and any ${\bf u}\in\Omega_\varphi$ we can form a
tiling of the real line by intervals.  For each letter $i\in\mathcal{A}$ we assign a closed
interval of length $l_i$.  We refer to the vector ${\bf L}=(l_1,\ldots , l_n)$ as the length
vector.  We then ``lay" copies of these closed intervals down on the real line (so that they do
not overlap on their interiors) according to the order prescribed by ${\bf u}=\ldots
u_{-1}.u_0u_1\ldots$, with the placement of the origin given by the decimal point.  There is a
natural translation acting $\Gamma^t$ on a tiling $T$ which forms new tilings by simply moving the
origin by a distance $t\in\reals$ to the right.  For a fixed ${\bf u}$, we form a compact metric
space by taking the completion of the space of the translates of ${\bf u}$.  Here the completion
is taken with respect to the metric which considers two tilings to be close if they agree on a
large neighborhood about the origin after a small translation.  If the substitution is primitive,
then this space is independent of ${\bf u}$, is minimal and uniquely ergodic.  We refer to the
$\reals$-action on $(\Omega_\varphi, \Gamma^t, {\bf L})$ as the {\it tiling dynamical system}.  
It is topologically conjugate to the suspension flow over the $\ints$-action $(\Omega_\varphi,
\sigma)$, with the height function equal to $l_i$ over the cylinder $i$.  The system
$(\Omega_\varphi, \Gamma^t, {\bf L}, \mu)$ is then said to have {\it pure point spectrum} if there
is a basis of $L^2(\Omega_\varphi, \Gamma^t, {\bf L}, \mu)$ consisting of eigenfunctions for the
$\reals$-action.

A primitive substitution is of {\it Pisot type} if all of its non-Perron-Frobenius eigenvalues are
strictly between $0$ and $1$ in magnitude.  Solomyak and Sirvent \cite{sirsol} have shown, for a
substitution of Pisot type, that if the $\reals$-action has pure discrete spectrum then so does
the $\ints$-action.  A recent result of Clark and Sadun \cite{clasad} implies that this can be
strengthened to an ``if and only if" statement.

In fact, the result of Clark and Sadun has more general applications than just the Pisot case
and we briefly describe some aspects of their results.  First, however, we should give a few
observations in order that we may put their results in context.

We denote as ${\bf L}_\lambda$ the left eigenvector associated to the Perron-Frobenius eigenvalue
$\lambda$ of the transition matrix $A_\varphi$.  The choice of ${\bf L}_\lambda$ as a length vector
causes any tiling which is (combinatorially) fixed under $\varphi$ to be {\it self-similar}.  
The similarity map is expansion by $\lambda$ and the tile-substitution is the geometric
version of $\phi$.  On the other hand, the length vector ${\bf L}_1=(1, 1,\ldots, 1)$ is a natural
choice since the group action of the time-$1$ return map for the tiling dynamical system is
conjugate to the shift homeomorphism of the substitution dynamical system.  Other length
changes, though producing tiling dynamical systems conjugate to the suspension of the shift,
will produce no such group action.  It is easy to see that the $\reals$-action of the tiling
dynamical system with constant length vector ${\bf L}_1$ has pure discrete spectrum if and only if
the $\ints$-action does.  Therefore, questions about the relation of the $\reals$-action with
length vector $L$ to the $\ints$-action can be rephrased into questions about the conjugacy
of the $\reals$-actions of the systems using ${\bf L}$ and ${\bf L}_1$, respectively.

Clark and Sadun \cite{clasad} give explicit conditions for the $\reals$-action associated to
one length vector to be conjugate (up to an overall rescaling) to that of another via a
homeomorphism which preserves the combinatorics (i.e. is homotopic to the identity).  Assume
the length vectors ${\bf L}$ and ${\bf L}'$ have been rescaled to agree in the
Perron-Frobenius direction.  Then, the conjugacy occurs if and only if ${\bf
L}A_\varphi^k-{\bf L}'A_\varphi^n\ra 0$ as $n\ra\infty$.  Note that every length vector can
be written as a linear combination of vectors living in the Perron-Frobenius eigenspace, the
small eigenspaces (those associated to the eigenvalues of magnitude less than one) and the
large eigenspaces (those associated to eigenvalues of magnitude greater than or equal to one,
though we don't include the P.F.-eigenspace here).  In particular, after the rescaling, ${\bf
L} - {\bf L}'$ must avoid the large eigenspaces for the conjugacy to occur.  If the
substitution is of Pisot-type, then all systems length changes produce conjugate systems and
we arrive at the result above.

We point out a few special kinds of substitutions we will be considering.  A substitution
$\varphi$ is said to be {\it irreducible} if the characteristic polynomial of its
transition matrix is irreducible.  In particular, ${\bf L}_\lambda$ of an irreducible
substitutions is such that none of its entries are rationally related to one another.  
Lastly, a substitution $\varphi$ is said to have {\it constant length} if the number of 
letters appearing in $\varphi (i)$ is the same for all $i\in\mathcal{A}$.

\section{Balanced Pair Algorithms}\label{balpairalgo}

We now describe the balanced pair algorithm in a variety of circumstances.  We begin with the
case that the substitution is irreducible.  This case was studied extensively in
\cite{sirsol} as an adaptation of the algorithm of \cite{liv2}.  The extension of this
algorithm to a specific class of reducible substitutions containing ``letter equivalences"
can then easily be described.  We end this section with description of how one extends this
procedure for generic reducible cases.

\subsection{Irreducible case:} A pair of allowable words $u$ and $v$ are called {\it balanced}
if the have the same population vector.  We write $\left|{u\atop v}\right|$ if $u$ and $v$
are balanced.  Note that if $\left|{u\atop v}\right|$ is balanced, then so is
$\left|{\varphi(u)\atop\varphi(v)}\right|$.

Let the right infinite sequence ${\bf u}=u_0u_1\ldots$ be fixed under the substitution and let
$w$ be a non-empty prefix of ${\bf u}$.  Since $w$ appears in ${\bf u}$ with positive
frequency, ${\bf u}$ can be written as ${\bf u}=wX_1wX_2\ldots$.  We may then speak of splitting
$\left|{{\bf u}\atop\sigma^{|w|}{\bf u}}\right|$ into balanced pairs in the following way:
\[\left|{{\bf u}\atop\sigma^{|w|}{\bf u}}\right| = \left|{wX_1\atop 
X_1w}\right|\left|{wX_2\atop X_2w}\right|\cdots .\]  
Since appearances of $w$ in ${\bf u}$ are bounded, there are only finitely many different
balanced pairs encountered in the process above.  We may further reduce each of these to form
a finite set of irreducible balanced pairs which we will refer to as $I_1(w)$.

We may now inductively define, for $n>1$:
\begin{align*}
I_n(w)=&\left\{\left|{u\atop v}\right| : \left|{u\atop v}\right| \text{ appears as an
irreducible balanced pair}\right. \\ 
&\left.\text{\ \ in the reduction of}\left|{\varphi(x)\atop \varphi(y)}\right| \text{, for 
some }\left|{x\atop y}\right|\in I_{n-1}\right\}.
\end{align*}

Let $I(w)=\bigcup_{n=1}^\infty I_n(w)$.  If $I(w)$ is finite, then we say that the balanced
pair algorithm associated to a prefix $w$, or bpa-$w$, terminates.  Below, we will state how
this algorithm is used to determine pure discrete spectrum, though our main interest here is
in determining precisely when the algorithm terminates.  We illustrate the computation of 
$I(w)$ for a simple example.

\begin{example}\end{example}
Consider the substitution given by:
\begin{center}
\begin{tabular}{ccl}
1 & $\ra$ & 112 \\
2 & $\ra$ & 12
\end{tabular}
\end{center}
Then ${\bf u}=11211212112\ldots$ and we take $w=1$.  Then it is easy to see that:
\[ I_1=\left\{\left|{1\atop 1}\right| , \left|{12\atop 21}\right|\right\}. \]
Now,
\[ \left|{1\atop 1}\right| \ra \left|{1\atop 1}\right|\left|{2\atop
2}\right| \text{ and }
\left|{12\atop 21}\right| \ra \left|{1\atop 1}\right|\left|{12\atop
21}\right|\left|{1\atop 1}\right|\left|{2\atop 2}\right| \]
so that $I_2=\left\{\left|{1\atop 1}\right| , \left|{12\atop 21}\right| , 
\left|{2\atop 2}\right|\right\}$ and further, $I_2 = I_3 =\ldots =I(w)$.  Thus, 
the algorithm terminates.

\subsection{Substitutions involving letter equivalences:} 
The balanced pair algorithm as originally described by Livshits \cite{liv1, liv2} has an
additional feature.  Consider two letters $i, j\in\mathcal{A}$ to be equivalent if $\varphi^n(i)$
and $\varphi^n(j)$ have the same number or symbols for all $n\in\nats$.  Then, a pair of words is
balanced if both contain the same number of symbols from each equivalence class.  The algorithm
runs in the usual way, starting with an initial list of balanced pairs, substituting and reducing.  
It stops if no new balanced pairs are produced.  In particular, in the constant length case all
letters are equivalent so all irreducible balanced pairs are just pairs of symbols.  Thus, the
algorithm always terminates in this case.  This way, the algorithm includes Dekking's criterion
\cite{dek} in the constant length case. For any irreducible case, no two symbols are equivalent
and thus the algorithm runs just as before.

\begin{example}(Substitution of constant length)\end{example} 
Consider the substitution given by: 
\begin{center} 
\begin{tabular}{ccl}
1 & $\ra$ & 112 \\
2 & $\ra$ & 122
\end{tabular}
\end{center}
which is a constant length substitution.  If we ignore the equivalence relation of Livshits,
then the troubling balanced pair is $\left|{12\atop 21}\right|$.  Iterating this pair, we
see:

\[\left|{12\atop 21}\right| \ra \left|{1\atop 1}\right|\left|{1212\atop 2211}\right|\left|{2\atop 
2}\right|\]
\[\left|{1212\atop 2211}\right| \ra \left|{1\atop 1}\right|\left|{1212211212\atop 
2212211211}\right|\left|{2\atop 2}\right|.\]
In particular, this process generates new balanced pairs of the form $\left|{1w2\atop 
2w1}\right|$ for longer and longer words $w$.

Once we take the equivalence relation into account, however, all letters are equivalent and
therefore the above process terminates.

\begin{example}\label{exnoncon}(A non-constant length substitution)\end{example}
Consider the substitution given by:
\begin{center}
\begin{tabular}{ccl}
1 & $\ra$ & 31 \\
2 & $\ra$ & 412 \\
3 & $\ra$ & 312 \\
4 & $\ra$ & 412
\end{tabular}
\end{center}
If we were to again ignore the equivalence relation of Livshits, then there is a potential 
problem with balanced pairs that ``match up" the letter $3$ with $4$.  To see this, note that 
the right eigenvectors of the transition matrix in some sense describe the frequencies in 
which letters appear.  Here, there is one large eigenvalue other than the Perron-Frobenius 
eigenvalue.  The eigenvalue is $1$ and it has right eigenvector $(0, 0, 1, -1)^T$.  We can 
thus view a balanced pair of the form $\left|{3\ldots\atop 4\ldots}\right|$ as initially 
having an abundance of a $3$ and a lack of a $4$ on top.  Since this corresponds to the right 
eigenvector of $1$, this difference persists under substitution.  Note that we say this is a 
``potential" problem as this association to a large eigenvalue in itself will not force the 
bpa-$w$ to not terminate.  In this case, however, we show that this does in fact occur.  
After shifting the fixed word ${\bf u}=312\ldots$ with $w=31$, we see the balanced pair 
$\left|{31412\atop 41231}\right|$.  Iterating this pair, we see:

\[\left|{31412\atop 41231}\right| \ra \left|{3123141231412\atop 4123141231231}\right|.\]
Thus, we generate new balanced pairs of the form $\left|{3w412\atop 4w231}\right|$ for longer and longer 
words $w$.  The equivalence relation tells us however that $2, 3$ and $4$ are actually equivalent letters 
and thus the above process terminates as $\left|{3\atop 4}\right|$ and $\left|{4\atop 
2}\right|\left|{12\atop 31}\right|$ are balanced.  Note that the left Perron-Frobenius 
eigenvector for this system is $(1, \lambda, \lambda, \lambda)$ so that $2, 3$ and $4$ all 
have the same length in the self-similar system.

\subsection{Reducible case:} The example above gives some hint as to how to describe the
balanced pair algorithm in the generic reducible case.  For a system whose length of tiles
have rational relations which persist under substitution, it may be the case that entire
words should be identified even when no individual letters are.  We therefore introduce an
equivalence relation which exploits these rational relations.  Consider two words $v$ and $w$
to be equivalent if: 
\[ {\bf L}({\bf p}(\varphi^n(v))-{\bf p}(\varphi^n(w)))=0, \forall n\in\nats , \] 
or equivalently, if 
\[ {\bf L}(A_\varphi^n{\bf p}(v) - A_\varphi^n{\bf p}(w))=0, \forall n\in\nats . \] 
We write $v\sim_{\bf L} w$ to emphasize the dependence of the equivalence relation on choice
of ${\bf L}$.  In the case in which we use the left Perron-Frobenius eigenvector ${\bf
L}_\lambda$ for our lengths, the equivalent pairs correspond simply to the {\it geometric
balanced pairs} of \cite{sirsol}.  In fact, we will show below that in this case, the
algorithm is equivalent to the overlap algorithm and for a general ${\bf L}$ is meant to
bridge the gap between the balanced pair and overlap algorithms in the reducible case.

The algorithm once again runs in the usual way, with a pair of words considered balanced if
they are equivalent under our relationship above.  We denote this algorithm by {\it bpa$(w,
L)$} and the set of balanced (or equivalent) pairs by $I(w, {\bf L})$ to emphasize that their
is now an additional dependence on ${\bf L}$.

\begin{remark}
Our equivalence condition above that ${\bf L}({\bf p}(\varphi^n(v))-{\bf p}(\varphi^n(w)))=0, 
\forall n\in\nats$, could be weakened to allow ${\bf L}({\bf p}(\varphi^n(v))-{\bf 
p}(\varphi^n(w)))\ra 0$ as $n\ra\infty$, which essentially allows inconsequential length 
changes in the small eigenspaces mentioned above.
\end{remark}

\begin{example}\end{example}
Consider the substitution given by:
\begin{center}
\begin{tabular}{ccl}
1 & $\ra$ & 112 \\
2 & $\ra$ & 2321 \\
3 & $\ra$ & 12.
\end{tabular}
\end{center}
The eigenvalues of the transition matrix are $\frac{3\pm\sqrt{13}}{2}$ and $1$ while
${\bf L}_\lambda=(1, \frac{-1 + \sqrt{13}}{2}, \frac{5-\sqrt{13}}{2})$.  The right eigenvector
associated to the eigenvalue $1$ is $(2, -1, -1)^T$.  This indicates that there is a
potential problem with balanced pairs which ``match up" the word $11$ with $23$.  Such a
matching occurs in the balanced pair $\left|{11223\atop 23211}\right|$, which was found by
performing the balanced pair algorithm after shifting the fixed word by two spaces.  
Iterating this balanced pair, we generate balanced pairs of the form $\left|{11w23\atop
23w11}\right|$ for longer and longer words $w$.

Fortunately, ${\bf L}_\lambda\cdot (2, -1, -1)=0$ so that $11$ and $23$ are equivalent words and 
hence $\left|{11\atop 23}\right|$ is balanced and the above algorithm will terminate.  Also 
note that ${\bf L}_1$ is also perpendicular to $(2, -1, -1)$ so that the algorithm will also 
terminate for what we consider to be the other interesting case.  This is immediate from the 
fact that these two systems are conjugate by \cite{clasad}.  On the other hand, ${\bf L}=(1, 
1, 2)$ for example is not conjugate to these and furthermore the algorithm does not terminate 
in this case.

\begin{proposition}\label{equivprop}
Let $\mathcal{E}(L)=\{ (u,v) : u\sim_{\bf L} v\}$.  Then, for any length vector ${\bf L}$, 
$\mathcal{E}(L)\subseteq\mathcal{E}({\bf L}_\lambda)$
\end{proposition}

\begin{proof}
Let $u\sim_{\bf L} v$.  Let ${\bf z}={\bf p}(u) - {\bf p}(v)$.  Since ${\bf L}$ has only
positive terms, it must have a component in the ${\bf L}_\lambda$ direction.  Let ${\bf
L}=a_\lambda {\bf L}_\lambda + \sum_{j=1}^lb_j B_j + \sum_{j=1}^sa_j A_j$, where $B_j$ is a
left eigenvector for eigenvalue $\beta_j\geq 1$ and $S_j$ is a left eigenvector for
eigenvalue $\alpha_j <1$ for each $j$.  Then, \[ 0= {\bf L}\cdot A^n{\bf z} =
(\lambda^na_\lambda {\bf L}_\lambda + \sum_{j=1}^l\beta_j^nb_j B_j +
\sum_{j=1}^s\alpha_j^na_j A_j)\cdot {\bf z}. \] If ${\bf L}_\lambda\cdot {\bf z}\not= 0$,
then the above implies $\lambda^n\approx\sum_{j=1}^lC_j\beta_j^n$ for large enough $n$ and
some constants $C_j$.  This is a contradiction since the Perron-Frobenius eigenvalue
dominates all others.
\end{proof}

\begin{corollary}\label{lambdaterm}
For any length vector ${\bf L}$, if the bpa$(w, {\bf L})$ terminates, then the bpa$(w, {\bf 
L}_\lambda)$ terminates.
\end{corollary}

Let $\mathcal{L}$ denote the span of the right eigenvectors with eigenvalues greater than or
equal to $1$ in magnitude but strictly less than the Perron-Frobenius eigenvalue.  
Similarly, let $\mathcal{S}$ denote the span of the right eigenvectors with eigenvalues
strictly less than $1$ in magnitude.  We say that $(u, v)$ lies in a vector space
$\mathcal{P}$ if $({\bf p}(u)-{\bf p}(v))\in\mathcal{P}$.  

\begin{remark} 
Then, the set of equivalence words lies entirely in the spaces since they form precisely what
is perpendicular to ${\bf L}_\lambda$.  In our examples above, choosing an ${\bf L}$ which
missed identifying equivalent pairs $\mathcal{L}$ would have led to the bpa$(w, {\bf L})$ not
terminating.  Contrast this with Example \ref{exnoncon} in which choosing an ${\bf L}$ which
neglects to identify $2$ and $4$ will have no effect on whether the algorithm terminates.  
The vector associated to the pairing of $2$ and $4$, namely $(0, 1, 0, -1)$ lives in the
zero-eigenspace and thus differences in frequencies should be quickly dispelled.  Our
contention here is that generally the equivalence relations which live in $\mathcal{S}$
should not affect the algorithm whereas those in $\mathcal{L}$ affect it greatly.
\end{remark}

\begin{corollary}
$u\sim_{\bf L} v$ implies ${\bf p}(u) - {\bf p}(v)\in\mathcal{L}\oplus\mathcal{S}$.  The 
converse 
is true if ${\bf L}={\bf L}_\lambda$.
\end{corollary}

A balanced (equivalent) pair $\left|{i\atop i}\right|$, for $i\in\mathcal{A}$, is called a
{\it coincidence}.  We say that a balanced pair $\left|{u\atop v}\right|$ {\it leads to a
coincidence} if there exists $m$ such that the reduction of
$\left|{\varphi^n(u)\atop\varphi^n(v)}\right|$ contains a coincidence.  Notice that
coincidences lead to coincidences since $\left|{\varphi(i)\atop\varphi(i)}\right|$ has
nothing but coincidences in its reduction.

\begin{theorem}\label{thmpds}
Let $\varphi$ be a primitive substitution such that ${\bf u}=u_0u_1\ldots$ is a right 
infinite fixed word and let ${\bf L}$ be a length vector.

\begin{enumerate}

\item[(a)] If for some prefix $w$ the bpa$(w, {\bf L})$ terminates and every equivalent pair
in $I(w, {\bf L})$ leads to a coincidence, then $(\Omega_\varphi, \Gamma^t, {\bf L}_\lambda)$
has pure discrete spectrum.

\item[(b)] If the bpa$(w, {\bf L})$ terminates for some prefix $w=u_0\ldots u_m$ such that
$u_{m+1}=u_0$, and $(\Omega_\varphi, \Gamma^t, {\bf L}_\lambda)$ has pure discrete spectrum, then
every balanced pair in $I(w, {\bf L})$ leads to a coincidence.


\end{enumerate}

\end{theorem}

Before beginning the proof of Theorem \ref{thmpds}, we make the following observations
regarding the densities of coincident pairs. Let ${\bf u}=u_0u_1\ldots$ be a fixed
right-sided sequence. Let $z$ be a prefix of ${\bf u}$ and ${\bf u}-{\bf L}\cdot{\bf
p}(z)=v_0v_1\ldots$. Let $D(z)=\{ u_i : u_i=v_j$ some $j$ with $u_0\ldots u_{i-1}\sim_{\bf L}
v_0\ldots v_{j-1}\}$ Suppose we define a density function,
$dens_{\bf L}(D(z))=\lim_{k\ra\infty}\frac{{\bf L}\cdot{\bf p}(D(z)\cap u_0\ldots u_k)}{{\bf
L}\cdot{\bf p}(u_0\ldots u_k)}$, if the limit exists.  The existence of this limit follows
from the unique ergodicity of $(\Omega_\varphi, \Gamma^t, {\bf L})$.  Notice that for ${\bf
L}={\bf L}_1$, this definition of density agrees with that used in \cite{sirsol} for the
(irreducible) balanced pair algorithm and for ${\bf L}={\bf L}_\lambda$ it agrees with that
of \cite{sol} for the overlap algorithm.  Now, by Proposition \ref{equivprop}, coincident
pairs for $\mathcal{E}({\bf L})$ are also coincident pairs for $\mathcal{E}(L_\lambda)$.  

\begin{remark}
The proof of this theorem differs from that of the irreducible case only in the way in which
we define density and the set of irreducible balanced pairs.  We therefore only include a
sketch of the proof that follows closely a sketch provided in \cite{sirsol} for the
irreducible case.  The full details of that case has been worked out in \cite{hol}.  A 
theorem of this sort was proved by \cite{liv2}, though coincidences and the balanced pair 
algorithm go back to \cite{dek} and \cite{mic}, respectively.  Part $(a)$ is largely 
contained in \cite{quef}.
\end{remark}

\begin{proof} (of Theorem \ref{thmpds}) 

Let $w$ be a prefix of the fixed word $\bf{u}$.  Denote by $D(z)$ the density defined above.  
We will be interested in $dens_{\bf L}(D(\varphi^l(w)))$ as $l\ra\infty$.  Let
\[ {\bf u}^{(l)}=\mu_1^{(l)}\mu_2^{(l)}... \] 
be the reduction of $\left|{\bf u}\atop\sigma^{p_l}{\bf u}\right|$ into irreducible 
equivalent pairs. For an equivalent pair $\beta =\left|{u\atop v}\right|$ let $|\beta|={\bf 
L}\cdot {\bf p}(u)$ and $\delta (\beta )=\{ u_i\in u : u_i=v_j$ some $j$ with $u_0\ldots 
u_{i-1}\sim_{\bf L} v_0\ldots v_{j-1}\}$.  Then 
\[ dens_{\bf L}(D(\varphi^l(w))) = \lim_{N\ra\infty}\frac{\sum_{j=1}^N{\bf L}\cdot{\bf 
p}(\delta (\mu_j^{(l)}))}{\sum_{j=1}^N{\bf L}\cdot{\bf p}(\mu_j^{(l)})}. \] 
Now consider the substitution $\hat\varphi$ on the set of irreducible balanced pairs $I(w,
{\bf L})$.  By definition clearly ${\bf u}^{(l)}\in I(w, {\bf L})^\nats$ and ${\bf
u}^{(l)}=\hat\varphi^l({\bf u}^{(0)})$, where ${\bf u}^{(0)}$ is the reduction of
$\left|{{\bf u}\atop\sigma^{|w|}{\bf u}}\right|$ into irreducible equivalent pairs.

There is a directed graph $\mathcal{G}(\hat\varphi)$ associated with the substitution
$\hat\varphi$.  Its vertices are labelled by the members of $I(w, {\bf L})$, and for every
vertex $\beta$ there are directed edges from $\beta$ into of the letters of $\hat\varphi
(\beta)$ with multiplicities.

(a) Let $\beta\in I(w, {\bf L})$ be an irreducible equivalent pair which is not a
coincidence. By assumption there is a path in the graph $\mathcal{G}(\hat\varphi)$ leading
from $\beta$ to a coincidence.  Since all the edges from coincidences lead to coincidences, a
standard argument shows that the frequency of the symbol $\beta$ in ${\bf
u}^{(l)}=\hat\varphi^l({\bf u}^{(0)})$ goes to zero geometrically fast as $l\ra\infty$.  
Since $I(w, {\bf L})$ is finite and $1-\delta (\beta)>0$ if and only if $\beta$ is a
non-coincidence, it follows that

\[ 1-dens(D_{p_l})\leq const\cdot\gamma^l \]
for some $\gamma\in (0, 1)$.  This implies that $\bf u$ is mean-almost periodic so we can
conclude (similar to \cite{quef}), VI.25) that $\varphi$ has pure discrete spectrum.  We note
that the argument only relied on $I(w, {\bf L})$ being finite and that the equivalence
relation does not create false coincidences.

(b)  Let $w=u_o...u_m$ be such that $u_{m+1}=u_0$ and $p_l={\bf L}\cdot{\bf
p}(\varphi^l(w))$.  It follows from \cite{host} (see also \cite{sol}, Theorem 4.3) that
$\lim_{l\ra\infty}e^{2\pi i\lambda p_l}=1$ for any eigenvalue $e^{2\pi i\lambda}$ of the
dynamical system $(\Omega_\varphi , \Gamma^t, {\bf L}, \mu)$.  If the spectrum is pure
discrete, then the eigenfunctions span a dense subset of $L^2(\Omega_\varphi )$ so that
$\lim_{l\ra\infty}||U_\varphi^{p_l}f-f||_2=0$ for every $f\in L^2(\Omega_\varphi )$.  Taking
$f$ to be the characteristic function of the cylinder set corresponding to $i\in\mathcal A$
with heights prescribed by ${\bf L}$, we obtain (see \cite{quef}, just after Lemma VI.26)  
$\lim_{l\ra\infty}dens_{\bf L}(D(\varphi^l(w)))=1$.  On the other hand, suppose that there is
an irreducible equivalent pair in $I(w, {\bf L})$ which does not lead to a coincidence.  
Then there exists an irreducible component $\mathcal{G}_0$ of the graph
$\mathcal{G}(\hat\varphi)$ which contains no coincidences.  There exists $l_0$ such that for
every $l\geq l_0$ elements of the component $\mathcal{G}_0$ occur in ${\bf u}^{(l)}$ with
positive frequency.  (Note that different elements of $\mathcal{G}_0$ may occur for different
$l$).  Further, it can be shown that this frequency is bounded away from zero as
$l\ra\infty$.  Since $1-\delta(\beta)>0$ for all $\beta\in\mathcal{G}_0$, it follows that
$dens_{\bf L}(D(\varphi^l(w)))\not\ra 1$, which is a contradiction.

\end{proof}

We also note the following relationship between the balanced pair algorithm and the overlap 
algorithm in the case that the tiling space is self-similar.

\begin{proposition}\label{propoverlap}
Let $\varphi$ be a primitive substitution such that $\varphi(1)$ begins with $1$.
Then the overlap-algorithm associated to $x=x(w)$ terminates with half-coincidences if and 
only if bpa$(w,{\bf L}_\lambda)$ terminates.
\end{proposition}

\begin{proof}

Assume the bpa$(w,{\bf L}_\lambda)$ terminates.  Then the distance between half-coincidences
(endpoints of our equivalent pairs) arising from looking at $(T, T -\lambda^nx)$ is bounded,
where $x=x(w)={\bf L}\cdot{\bf p}(w)$.  Theorem 5.6 of \cite{sirsol} implies that the overlap
algorithm associated to $x(w)$ terminates with half-coincidences.

Assume the overlap-algorithm associated to $x={\bf L}\cdot{\bf p}(w)$ terminates with
half-coincidences.  Again by Theorem 5.6 of \cite{sirsol}, the distance between
half-coincidences arising from $(T, T -\lambda^nx)$ is bounded.  Hence bpa$(w, {\bf L}_\lambda)$
terminates.


\end{proof}

\begin{example}A non-terminating example\end{example}
We now give an example which will not terminate even with the extended equivalence relations presented 
here.  This example is a rewriting of the Morse-Thue systems in which we have forced 
coincidences.  We will use ${\bf L}_\lambda$ as our length vector so that we are conjugate to the 
original Morse-Thue system and hence do not have pure discrete spectrum.  This example 
therefore cannot terminate for any version of the balanced pair algorithm.

The substitution is given by: 
\begin{center}
\begin{tabular}{ccl}
1 & $\ra$ & 1234 \\
2 & $\ra$ & 124 \\
3 & $\ra$ & 13234 \\
4 & $\ra$ & 1324
\end{tabular}
\end{center}

Here, the P.F. eigenvalue is $4$ with left eigenvector $(3, 2, 4, 3)$.  (Note that this vector
also corresponds to the $\ints$-action of the original Morse-Thue system.)  The other important
eigenvalue (the remaining two are zero) is $1$ with right eigenvector $(1, 1, -2, 1)^T$.  Thus,
$abd$ is equivalent to $cc$, however these words do not cluster close enough to each other to aid
in terminating the balanced pair algorithm.  Notice also the system generated by the length vector
$(1, 1, 1, 1)$ is not conjugate to the system generated by ${\bf L}_\lambda$.  But by Corollary
\ref{lambdaterm}, the balanced pair algorithm will also not terminate for this length vector.  In
fact, using techniques from \cite{clasad} to directly compute its spectrum, one can see that it
will not have pure discrete spectrum either.

\section{Open problems and conjectures}\label{opensec}

The motivation for studying the effects of the new equivalence relation on the balanced pair
algorithm was not only to understand reducible substitutions, but was mainly to aid in
analyzing substitutions which have been collared or rewritten. We use the example below as a
test case for a general Pisot substitution.  Beginning with a Pisot substitution, we will
rewrite it into an equivalent substitution that always begins with the same letter.  In this
way we force every balanced pair to lead to a coincidence.  It remains to show that the new
substitution will terminate in order to show that all Pisot substitutions have pure discrete
spectrum.  A difficulty arises in that rewriting increases the size of the alphabet and can
therefore add additional eigenvalues of $0$ and $\pm 1$.  The zeros do not concern us, but
the roots of unity might.  Considering the $\ints$-action on our original substitution by
changing tiles to unit lengths produces an integer vector in the rewritten substitution which
will not have a component in the eigenspace of the roots of unity.  Further, these
eigenspaces generate equivalent words so that it is our hope that the equivalence classes
will always force the algorithm to terminate.

\begin{example}A rewritten Pisot substitution\end{example}
Consider the substitution $\varphi$ given by $a\ra abb$ and $b\ra ba$.  Using the rewriting
procedure of \cite{bardia1}, we square $\varphi$ and generate a substitution on $1=abbb, 2=ab,
3=aabbb$ and $4=aabb$ by:
\begin{center}
\begin{tabular}{ccl}
1 & $\ra$ & 122334 \\
2 & $\ra$ & 1224 \\
3 & $\ra$ & 12322334 \\
4 & $\ra$ & 1232234
\end{tabular}
\end{center} 
The eigenvalues of the transition matrix of this substitution are $0, 1\pm\frac{\sqrt{5}}{2}$
and $1$.  Since the original substitution is Pisot, it is insensitive to length changes in
$a$ and $b$.  The new substitution will also be insensitive to any length change which is
consistent with length changes in $a$ and $b$.  For example, setting $a=b=1$ generates the
length vector ${\bf L}=(3, 2, 5, 4)$ and will in fact miss the eigenspace associated to $1$.  
(Note, however, that ${\bf L}_1=(1, 1, 1, 1)$ does not.)  The right eigenvalue of $1$,
$v_1=(1, 1, -2, 1)^T$, is perpendicular to ${\bf L}$ so that $abd\sim_{\bf L} cc$.  Because
${\bf L}_1$ is not perpendicular to $v_1$, this makes the bpa$(w, {\bf L})$ difficult to run.  
Some tedious calculations reveal that bpa$(w, {\bf L})$, where one must keep an eye out for
the equivalence relations, will terminate.  This example produces $30$ different irreducible
balanced pairs, the longest one of which contains a word of length $11$.  We suspect that
bpa$(w, {\bf L}_1)$ will not terminate, but this has also not yet been shown.

Generally, if the balanced pair algorithm terminates for a Pisot substitution, \cite{sirsol} 
gives us that the distance between half-coincidences is bounded.  For the rewritten 
substitution, the new tiles are compositions of smaller old tiles and these half-coincidences 
may occur ``internally."  We suspect this cannot happen and that the algorithm must terminate 
for the new substitution as well.  More precisely we have:

\begin{conjecture}
Let $\varphi$ be a Pisot substitution on $n$-letters.  Let $\tilde\varphi$ be a rewriting of
$\varphi$ so that or some letters $b, e$ in the rewritten alphabet, $\tilde\varphi(i)=b\ldots e$
for all $i$.  Assume that the balanced pair algorithm for the original Pisot substitution
$\varphi$ terminates.  Then the balanced pair algorithm for the rewritten substitution
$\tilde\varphi$ also terminates for length vector $L_\lambda$.
\end{conjecture}

An immediate corollary of this would be that if the balanced pair algorithm of a Pisot type 
substitution terminates, then it must do so with coincidences.

{\bf Acknowledgements:} The author would like to thank Lorenzo Sadun for his numerous
contributions to this work and helpful suggestions.  Additionally, we thank Boris Solomyak for
some especially illuminating correspondences.

\end{document}